\documentclass[12pt]{amsart}
\usepackage{amssymb}

\newtheorem{theorem}{Theorem}

\newtheorem{proposition}{Proposition}
\newtheorem{corollary}{Corollary}
\newtheorem{remark}{Remark}

\newtheorem{conjecture}{Conjecture}
\newtheorem{lemma}{Lemma}

\newcommand{\Z}{{\mathbb Z}}
\newcommand{\Q}{{\mathbb Q}}
\newcommand{\R}{{\mathbb R}}
\newcommand{\C}{{\mathbb C}}

\newcommand{\RP}{{\mathbb RP}}
\newcommand{\CP}{{\mathbb CP}}

\begin{document}

\title{Integrable bodies in odd-dimensional spaces}

\author{V.A.~Vassiliev}
\address{Steklov Mathematical Institute of Russian Academy of Sciences \ \ and \ \ National Research University Higher School of Economics}
\thanks{Research supported by the Russian Science Foundation grant, project
16-11-10316}
 \email{vva@mi-ras.ru}



\begin{abstract}
V. Arnold's problem 1987-14 from \cite{A} asks whether there exist bodies with smooth boundaries in $\R^N$ (other than the ellipsoids in odd-dimensional spaces) for which the volume of the segment cut by any hyperplane from the body is an algebraic function of the hyperplane. We present a series of very realistic candidates for the role of such bodies, and prove that the corresponding volume functions are at least algebroid, in particular their analytic continuations are finitely valued; to prove their algebraicity it remains to check the condition of finite growth. 
\end{abstract}

\maketitle

\section{Introduction}

Any compact body with regular boundary  in ${\mathbb R}^N$  defines a two-valued function on the space of affine hyperplanes: its values are the volumes of two parts into which the hyperplanes cut the body.

A body in $\R^N$ (and also the hypersurface bounding it) is called {\it algebraically integrable} if this function is algebraic, i.e. there is a non-zero polynomial $F$ in $N+2$ variables such that $F(V,a_1, \dots, a_N,b)=0$ any time when $V$ equals the volume cut (to any side) from this body by the hyperplane defined by the equation \ $a_1x_1 + \cdots + a_Nx_N+b=0$.

Archimedes \cite{arch} has proved that spheres in $\R^3$ are algebraically integrable; today it is easy to check this property also for all ellipsoids in odd-dimensional spaces. Newton (\cite{Newton}, Lemma XXVIII) has proved that bounded convex bodies with smooth boundaries in $\R^2$ never are algebraically integrable. V.~Arnold asked (see \cite{A}, Problems 1987-14, 1988-13, and 1990-27, and also \cite{Hooke}) whether there exist algebraically integrable bodies with smooth boundaries other than collections of ellipsoids in odd-dimensional spaces. It was proved in \cite{V14} that there are no such bodies in even-dimensional spaces. The main obstruction to the integrability is the monodromy action on homology groups related with hyperplane sections of our hypersurfaces: this action controls the ramification of the analytic continuation of volume functions, and usually provides infinitely many values of such a continuation at one and the same hyperplane. This obstruction implies even a more strong assertion: the volume function of a body with smooth boundary in $\R^N$, $N$ even, cannot be algebroid. It provides also many restrictions on the geometry of algebroidally integrable bodies for odd $N$ proving that such bodies  other than quadrics are very rare, see \cite{APLT}.

We present below a new series of bodies in all spaces of odd dimensions greater than three, for which this obstruction vanishes: it are tubular neighborhoods of standard even-dimensional spheres (and also any bodies affine equivalent to them). It is extremely plausible that these bodies are algebraically integrable: we prove that the analytic continuations of their volume functions are algebroid (in particular, finitely valued), so it remains only to prove that they have  finite  growth at their singular points.

\subsection{Main result}
Let $n$ be odd, $m$ be even, and $\varepsilon \in (0,1)$. Consider the Euclidean space $\R^{n+m} \equiv \R^n \oplus \R^m$ with orthogonal coordinates $x_1, \dots, x_n$ in $\R^n$ and $y_1, \dots, y_m$ in $\R^m$, and the hypersurface in it defined by equation
\begin{equation}
\left( \sqrt{x_1^2 + \dots + x_n^2}-1\right)^2 +y_1^2 + \dots + y_m^2 = \varepsilon^2 \ .
\label{meq}
\end{equation}

Denote by ${\mathcal P}$ (respectively, by ${\mathcal P}^{\C}$) the space of all affine hyperplanes in $\R^{n+m}$ (respectively, in $\C^{n+m}$), and by ${\mathcal Reg}$ the subset in ${\mathcal P}$ consisting of hyperplanes transversal to the hypersurface (\ref{meq}). Identifying any real hyperplane with its complexification, we consider ${\mathcal P}$ as a subset of ${\mathcal P}^{\C}$.

\begin{theorem}
\label{mthm}
For any \ $\varepsilon \in (0,1)$, there is a four-valued analytic function on the space ${\mathcal P}^{\C}$ such that one of volumes cut by hyperplanes from the body bounded by hypersurface $($\ref{meq}$)$ in $\R^{n+m}$  coincides in any point $X \in {\mathcal P}$ with a sum of at most two $($maybe zero$)$ values of this function. The monodromy group of this function is the Klein four-group generated by permutations $\binom{1234}{2143}$ and $\binom{1234}{4321}$ of values.
\end{theorem}

\begin{conjecture}
\label{conn1}
This four-valued analytic function is in fact algebraic.
\end{conjecture}

\begin{remark} \rm
If $n=1$ or $m=0$ then (\ref{meq}) defines a reducible hypersurface, and hence does not fit the Arnold's question.
\end{remark}

\begin{remark} \rm
Hypersurfaces defined by similar equations
\begin{equation}
\left( x_1^2 + \dots + x_n^2-1\right)^2 +y_1^2 + \dots + y_m^2 = \varepsilon^2
\label{loca}
\end{equation}
 with odd $n$ and arbitrary $m$ were considered in \cite{ne} as boundaries of {\em locally} algebraically integrable bodies. The hypersurface (\ref{meq}) has a smaller monodromy group of homology classes of hyperplane sections (which in fact controls the ramification of volume function) than (\ref{loca}), because it has singular points in the imaginary domain; these singularities ``eat''  parabolic points of the hypersurface, which cause a large part of the corresponding monodromy group.
\end{remark}

\begin{remark} \rm
It is important that the second fundamental form of the hypersurface (\ref{meq}) in $\R^{n+m}$ has only even inertia indices at all points where it is non-degenerate, cf. \S 3 in \cite{VA}.
\end{remark}

\begin{remark} \rm
Are there any other examples of this sort? Some natural candidates are provided by the iteration of the previous construction: we take a flag of odd-dimensional spaces $\R^k \subset \R^n \subset \R^N$ and a tubular neighbourhood  in $\R^N$ of the boundary of a tubular neighbourhood in $\R^n$ of the standard sphere $S^{k-1} \subset \R^k$.
\end{remark}

\begin{remark} \rm
Even if Conjecture \ref{conn1} will be proved, 
one of Arnold's problem of this series will remain unsolved: it asks whether there are extra {\em convex} algebraically integrable bodies in $\R^N$, see \cite{A}, problem 1990-27. Also, we do not have good examples for $N=3$.
\end{remark}

\begin{remark} \rm
See \cite{Agr}, \cite{KMY} for some related problems and results.
\end{remark}

An obvious problem (whose solution would imply Conjecture \ref{conn1})  is to integrate explicitly in radicals the volume function of the body (\ref{meq}).

\subsection{A detalization of Theorem \ref{mthm}}
\label{mlmm}

Denote by $W$ the body in $\R^{n+m}$ bounded by hypersurface (\ref{meq}); denote its volume by $C_0$, and the corresponding two-valued volume function on ${\mathcal P}$ by $V_W$. Denote by ${\mathcal P_2}$ (respectively, by ${\mathcal P}_2^{\C}$) the two-dimensional subspace in ${\mathcal P}$ (respectively, in ${\mathcal P}^{\C}$) consisting of all hyperplanes defined by equations of the form
\begin{equation}
px_1+qy_1+r =0
\label{hhh}
\end{equation}
with real (respectively, complex) coefficients $p, q, $ and $r$; also denote by ${\mathcal Reg}_2$ the space ${\mathcal P}_2 \cap {\mathcal Reg}$. \medskip

The group $O(n) \times O(m)$ of independent rotations in $\R^n$ and $\R^m$ acts on the space ${\mathcal P}$. The volume function $V_W$ is constant on the orbits of this action. Any such orbit contains points of the space ${\mathcal P}_2$.
In particular, any hyperplane defined by equation
\begin{equation}
\alpha_1 x_1 + \dots + \alpha_n x_n + \gamma_1 y_1 + \dots + \gamma_m y_m = \beta \ ,
\label{eqhy}
\end{equation}
where not all coefficients $\alpha_j$ are equal to $0,$ can be reduced by this group to a hyperplane with equation
\begin{equation}x_1 = ay_1+c\  ,
\label{hhhr}
\end{equation}
where
\begin{equation}
a = \frac{\sqrt{\gamma_1^2+ \dots +\gamma_m^2}}{\sqrt{\alpha_1^2+ \dots +\alpha_n^2}}, \quad c= \frac{|\beta|}{\sqrt{\alpha_1^2+ \dots + \alpha_n^2}} \ ,
\label{nof}
\end{equation}
cf. \cite{ne}.

The space ${\mathcal Reg}_2$ consists of five connected components. Indeed, the coordinate plane $\R^2 \subset \R^{n+m}$ defined by conditions $x_2 =\dots =x_n = y_2 = \dots = y_m=0$ intersects  the boundary $\partial W$ of the body $W$ along two circles \ $(x_1-1)^2 + y_1^2 =\varepsilon^2$ \ and \ $(x_1+1)^2 + y_1^2 =\varepsilon^2,$ see Fig. \ref{slice}. For any hyperplane $X \in {\mathcal Reg}_2$, the line $X \cap \R^2$ can

(1) not intersect these circles leaving them to one side of it,

($2_l$), ($2_r$) intersect transversally only the left-hand (respectively, right-hand) circle,

(3) not intersect these circles but separate them, or

(4) intersect both of them.
\medskip

\unitlength 1.20mm
\linethickness{0.4pt}
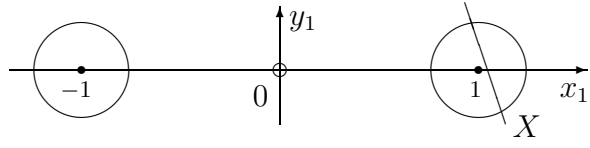
\begin{figure}
\begin{picture}(64,13)
\put(0,6){\vector(1,0){64}}
\put(61,3){$x_1$}
\put(30,6){\circle{1.5}}
\put(27,2){$0$}
\put(5.7,3){{\scriptsize $-1$}}
\put(51,3){{\scriptsize $1$}}
\put(8,6){\circle*{1}}
\put(52,6){\circle*{1}}
\put(8,6){\circle{10}}
\put(52,6){\circle{10}}
\put(55,0){\line(-1,3){4.5}}
\put(30,0){\vector(0,1){13}}
\put(31,11){$y_1$}
\put(55.8,-1.5){$X$}
\end{picture}
\caption{Planar slice of $W$}
\label{slice}
\end{figure}

Correspondingly, the space ${\mathcal Reg}$ splits into {\em four} components filled in by $(O(n)\times O(m))$-orbits of points from these components of ${\mathcal Reg}_2$; the components ($2_l$) and ($2_r$) of ${\mathcal Reg}_2$ define one and the same component (2) of ${\mathcal Reg}$. Then Theorem \ref{mthm} has the following detalization.

\begin{theorem}
\label{redmt}
There is a four-valued analytic function $\Psi$ on the space ${\mathcal P}^{\C}$ such that 

1$)$ its restriction to ${\mathcal P}$ is $(O(n) \times O(m))$-invariant,

2$)$ the sum of its four values is everywhere equal to $C_0$;

3$)$ any value of the volume function $V_W$ on any hyperplane $X \in {\mathcal Reg}$ coincides with

-- a value of this function $\Psi$, or $C_0$ less such a value, if $X$ belongs to the domain  $(2)$ of ${\mathcal Reg}$,

-- the sum of some two values of $\Psi$ if $X$ belongs to one of domains $(3)$ or $(4)$.

4$)$ in the latter case of domains $(3)$ or $(4)$, two values of the function $V_W$ can be obtained one from the other by analytic continuation along a path lying inside this domain. On the contrary, two leaves of the function $V_W$ in the domain $($2$)$ of ${\mathcal Reg}$ are not analytic continuations of one another even through the complex space ${\mathcal Reg}^{\C}$;

5$)$ the monodromy group of this analytic function is isomorphic to $\Z_2 \oplus \Z_2$ and is generated by permutations $\binom{1234}{2143}$ and $\binom{1234}{4321}$ of values;

6$)$ the restriction of this function $\Psi$ to the subspace ${\mathcal P}_2^{\C}$ splits into the product of two two-valued algebroid functions, one of which has the form 
\begin{equation}
\Phi(V,a,c) \equiv V^2 - S(a,c) V + P(a,c) 
\label{phi}
\end{equation}
 in coordinates $a, c$ $($see $(\ref{hhhr}))$, 
where $S$ and $P$ are single-valued analytic functions invariant under reflection $a \mapsto -a$, and the other one is equal to $\Phi(V,-a,-c)$.
\end{theorem}

\subsection{Scheme of the proof of Theorem \ref{redmt} (cf. \cite{V14}).}
\label{scheme}

Denote by $A$ the hypersurface in $\C^{n+m}$ defined by equation (\ref{meq}) or, equivalently, by equation
\begin{equation}\left(x_1^2 +\dots +x^2_n + y^2_1+ \dots + y_m^2+ (1-\varepsilon^2)\right)^2 - 4(x_1^2+\dots + x_n^2) =0 \ .
\label{sme}
\end{equation}
Denote  by ${\mathcal Reg}^{\C}$ the subset in ${\mathcal P}^{\C}$ consisting of all complex hyperplanes in $\C^{n+m}$ generic with respect to $A$ (that is, of hyperplanes whose closures in the compactification $\CP^{n+m}$ of $\C^{n+m}$ are transversal to the stratified variety consisting of the hypersurface $A$ and the ``infinitely distant'' plane $\CP^{n+m} \setminus \C^{n+m}$). Denote by ${\mathcal Reg}^{\C}_2$ the intersection ${\mathcal Reg}^{\C} \cap {\mathcal P}_2^{\C}$.

By Thom isotopy lemma (see e.g. \cite{GM}), groups $H_{n+m}(\C^{n+m}, A \cup X)$ are isomorphic to one another for all $X \in {\mathcal Reg}^{\C}$; moreover, they form a covering over ${\mathcal Reg}^{\C}$ with canonical flat connection (Gauss-Manin connection, see e.g. \cite{AVG}, \cite{APLT}). This connection naturally defines the monodromy action $ \pi_1({\mathcal Reg}^{\C},X) \to \mbox{Aut} ( H_*(\C^{n+m}, A \cup X))$. Integrals of the volume form
\begin{equation}
\label{volform}
dx_1 \wedge \dots \wedge dx_n \wedge dy_1 \wedge \dots \wedge dy_m
\end{equation}
along the elements of all these homology groups are well-defined and provide linear functions $H_{n+m}(\C^{n+m}, A \cup X; \C) \to \C$.

Let $X_0$ be a distinguished point in ${\mathcal Reg}^{\C}$, and $\Xi$ an element of the group
\begin{equation}
H_{n+m}(\C^{n+m}, A \cup X_0) \ .
\label{rehom}
\end{equation}
 The pair $(X_0, \Xi)$ defines an analytic function on ${\mathcal Reg}^{\C}$: the value of its continuation along a path $l$ in ${\mathcal Reg}^{\C}$ connecting $X_0$ with some point $X'$ is equal to the integral of the form (\ref{volform}) along the cycle in $H_{n+m}(\C^{n+m}, A \cup X')$ obtained from $\Xi$ by the Gauss-Manin connection over our path $l$.

If $X_0 \in {\mathcal Reg}$ and $\Xi$ is the homology class of a part of the body $W$ cut from it by hyperplane $X_0$, then in a neighbourhood of the point $X_0$ in ${\mathcal P}$ this function coincides with the volume function; hence the analytic continuations of these two functions into the complex domain also coincide. The ramification of the analytic continuation of the volume function is thus controlled by the monodromy action of the group $\pi_1({\mathcal Reg}^{\C}, X_0)$ on the group  (\ref{rehom}).

In particular, suppose that this hyperplane $X_0$ belongs to ${\mathcal Reg}_2$ and is sufficiently close to the hyperplane given by equation $x_1=1$ (but is not equal to it, because the latter hyperplane is not generic ``at infinity''), so that the line $X_0 \cap \R^2$ intersects the right-hand circle of $\partial W$ only, see Fig. \ref{slice}. It is easy to see that the component of $W \setminus X_0$ containing the piece of $\R^2 \cap W$ placed to the right from $X_0$ in Fig. \ref{slice} is a {\em vanishing cycle} in the group (\ref{rehom}): it contracts to a point when we move the plane $X_0$ to the right parallel to itself  until the tangency with $\partial W$. We prove below that the orbit of this element of (\ref{rehom}) under the action of the group $\pi_1({\mathcal Reg}^{\C},X_0),$ consists of four elements, and this orbit splits into two two-element orbits of the action of  the group $\pi_1({\mathcal Reg}_2^{\C},X_0)$. Therefore the analytic continuation of the corresponding volume function from a neighbourhood of $X_0$ to entire ${\mathcal Reg}^{\C}$ (respectively,  ${\mathcal Reg}^{\C}_2$) is four-valued (respectively, two-valued). 

These two analytic continuations are exactly the functions $\Psi$ and $\Phi$ promised in Theorem \ref{redmt}.
The coefficients $S(a,c)$ and $P(a,c)$ in (\ref{phi}) are respectively the sum and the product of both values of the second analytic continuation at the point $X(a,c) \in  {\mathcal Reg}^{\C}_2$ defined by equation (\ref{hhhr}). By definition, these functions $S$ and $P$ are single-valued.

If a hyperplane $X$ with equation (\ref{hhhr}) intersects only the right-hand (respectively, left-hand) circle of $\partial W \cap \R^2$ then the volume of one of parts cut by it from the body $W$ is equal to one of roots of the polynomial $\Phi(\cdot, a,c)$ (respectively,  $\Phi(\cdot, -a,-c)$). If $X$ separates two circles of $\partial W \cap \R^2$ then any part of $W$ cut by $X$ is equal to the sum of both elements of the same orbit of the action of  $\pi_1({\mathcal Reg}_2^{\C},X_0)$ on the $\pi_1({\mathcal Reg}^{\C},X_0)$-orbit of the vanishing cycle, and hence its volume is equal to $S(a,c)$ or $C_0 - S(a,c)$. A rotation of $X$ inside this component of ${\mathcal Reg}$ permutes the half-spaces of $\R^{n+m}$ separated by $X$, and hence permutes also these two leaves of the function $V_W$. If $X$ intersects both these circles then any of these two parts is equal to the sum of two elements of the $\pi_1({\mathcal Reg}^{\C},X_0)$-orbit of our vanishing cycle which belong to different $\pi_1({\mathcal Reg}_2^{\C},X_0)$-orbits, so that the volume function at $X$ is equal to the sum of one root of $\Phi(\cdot,a,c)$ and one root of $\Phi(\cdot, -a,-c)$. The basic element of $\pi_1({\mathcal Reg}_2^{\C}, X_0)$ permutes these elements within any of two pairs, and hence also moves one of leaves of $V_W$ in the domain (4) of ${\mathcal Reg}$ to the other.
\medskip

All these facts will be proved in sections 2---4.

\subsection{On functions $S$ and $P$}

The functions $S(a,c)$ and $P(a,c)$ defined in this way are regular  in $\C^2$  outside the divisor given by equation
\begin{equation}
\label{sinsin}
a^2+1=0 \ ,
\end{equation}
as follows from the next proposition.

\begin{proposition}
Let $D$ be a compact subset in the space $\C^2$ of hyperplanes $($\ref{hhhr}$)$, and is separated from the divisor $($\ref{sinsin}$)$.
Then 

1$)$ there is a constant $C=C(D)$ such that for any hyperplane $X \in D$ the function $$\|(x,y)\|^2 \equiv |x_1|^2 + \dots +|x_n|^2 + |y_1|^2 + \dots + |y_m|^2$$ is regular on the stratified variety $A \cup X$ everywhere outside the ball $B_C \subset \C^{n+m}$ in which this function takes the values $\leq C$;

2$)$ the functions $|S|$ and $|P|$ are bounded on $D$.
\end{proposition}

\noindent
{\it Proof.}  Statement 1) follows immediately from equation (\ref{sme}). 

We can assume that $D$ is simply-connected: otherwise we cover it by finitely many simply-connected compact sets and prove our estimates for any of these sets separately. By Lemma \ref{redmt} and considerations of stratified Morse theory, all homology classes in $H_{n+m}(\C^{n+m},A \cup X)$, $X \in D$, which can be obtained by Gauss--Manin connection from our vanishing cycle, can be realized by compact cycles contained in $B_C$ and depending continuously on $X$. Therefore the integrals of the form (\ref{volform}) along all of them are uniformly bounded. \hfill $\Box$

\begin{conjecture} These functions $S$ and $P$ are rational; they
have poles of orders respectively $\leq (n+m-1)/2$  and $\leq (n+m)$ on the variety $($\ref{sinsin}$)$, and poles of orders $n+m$ and $2(n+m)$ on the divisor $\{c=\infty\}$, and no other singularities in ${\mathcal P}_2^{\C}$. In particular, the single-valued entire analytic functions $S(a,c)\times (a^2+1)^{(n+m-1)/2}$ and $P(a,c) \times (a^2+1)^{n+m}$  grow only polynomially in $\C^2$, hence they are polynomials, and function $($\ref{phi}$)$ is algebraic.
\end{conjecture}

This conjecture obviously implies Conjecture \ref{conn1}.

\section{Geometry and topology of variety (\ref{meq})}

\begin{lemma} \label{lem4}
1. The singular locus $\mbox{\rm sing} A$ of variety $A$ is distinguished by the system of equations
 \begin{equation} x_1 = \dots = x_n =0, \ y_1^2+ \dots + y_m^2 +1-\varepsilon^2 =0 \ ,
\label{singloc}
\end{equation}
in particular it is a smooth $(m-1)$-dimensional complex manifold.

2. Let $D$ be a small $2(n+1)$-dimensional open disc in $\C^{n+m}$ transversal to the manifold $\mbox{\rm sing}A$ $($e.g. a fiber of its tubular neighbourhood$)$. Then the variety $A \cap D$ is diffeomorphic to the zero set of a non-degenerate quadratic form in  $\C^{n+1}$, i.e. to the variety defined by equation
\begin{equation}
z_1^2 + \dots + z_{n+1}^2 =0 \ .
\label{quadrcone}
\end{equation}
\end{lemma}

\noindent {\it Proof} is immediate. \hfill $\Box$

\begin{corollary}
\label{cor2}
In conditions of Lemma \ref{lem4}, the group $H_i(A\cap D, A \cap \partial D)$ is equal to $\Z$ if $i$ is equal to $n+1$ or $n$, and is trivial for all other $i$. \hfill $\Box$
\end{corollary}

\begin{corollary}
\label{link}
The {\em complex link} $($see \cite{GM}$)$ of $\mbox{\rm sing} A$ is homotopy equivalent to $S^n$. \hfill $\Box$
\end{corollary}

\begin{lemma} \label{lem5}
$H_{n+m-1}(A) \simeq \Z^2 \ . $ This group is generated by

1$)$ the fundamental cycle of the manifold $\partial W \equiv A \cap \R^{n+m},$ and

2$)$ the set of points $(x_1, \dots, x_n, y_1, \dots, y_m) \in A$ such that all $x_j$ are real, $x_1^2 + \dots + x_n^2 \leq (1-\varepsilon)^2 ,$ the real parts of all $y_j$ are equal to $0$, and $y_1^2 + \dots + y_m^2 \geq \varepsilon^2-1$.

The group $H_{n+m-2}(A)$ is trivial.
\end{lemma}

\noindent
{\it Proof.}
Two described cycles are independent in $H_{n+m-1}(A)$ since they bound relative cycles in $H_{n+m}(\C^{n+m},A)$ whose volumes (i.e. the integrals of the form (\ref{volform}) along them) are not equal to zero and are incomparable for generic $\varepsilon$. Indeed, the volume of the tube $W$ tends to zero when $\varepsilon$ tends to 0, and the volume of the second cycle tends to a real non-zero  constant (whose sign is equal to $(-1)^{m/2}$), therefore the ratios of these two numbers run a continuum of values.

Denote by $X_1$ the (non-generic) hyperplane $\{x_1=0\}$ in $\C^{n+m}$. Consider the exact sequence of the pair $(A, A \cap X_1)$:
\begin{eqnarray}
\label{es}
\quad H_{n+m-1}(A \cap X_1) \to H_{n+m-1}(A) \to H_{n+m-1}(A,A \cap X_1) \to \\
\label{es1}
 H_{n+m-2}(A \cap X_1)  \to H_{n+m-2}(A) \to H_{n+m-2}(A, A \cap X_1)\ .
\end{eqnarray}
The left-hand group in (\ref{es}) is trivial as $A \cap X_1$ is a $(n+m-2)$-dimensional Stein space. The restriction of the real function $|x_1|$ to the non-singular variety $A \setminus X_1$ has exactly four Morse critical points with coordinates $ x_1= \pm1\pm \varepsilon,$ $x_2 = \dots = x_n=y_1=\dots=y_m=0;$ their Morse indices are equal to $n+m-1$. Therefore the concluding group $H_{n+m-2}(A,A \cap X_1)$ in (\ref{es1}) is trivial, and the concluding group $H_{n+m-1}(A,A \cap X_1)$ in (\ref{es}) is isomorphic to $ \Z^4$ and is generated by intersections of the above-described two cycles with half-spaces where $x_1\leq 0$ or $x_1 \geq 0$. The boundaries of two generators lying in the half-space $\{x_1 \geq 0\}$ are the cycles in $A \cap X_1$ which are independent in $H_{n+m-2}(A \cap X_1)$ by exactly the same reasons by which two cycles discussed in the present lemma are independent in $H_{n+m-1}(A)$. So the rank of the last arrow in (\ref{es}) is equal to 2, and its kernel (isomorphic to $H_{n+m-1}(A)$) is 2-dimensional, i.e. $H_{n+m-1}(A) \simeq \Z^2$. The space $A \cap X_1$ is analogous to $A$ in the space $X_1 \equiv \C^{n+m-1}$ instead of $\C^{n+m}$, which implies the equality of the group $H_{n+m-2}(A \cap X_1)$ to $\Z^2$. Two generators of this group belong to the image of the boundary map of our exact sequence, hence both arrows in (\ref{es1}) are trivial, and $H_{n+m-2}(A) \simeq 0$.  \hfill $\Box$ \medskip

Denote by $\AA$ the regular part $A \setminus \mbox{sing}A$ of variety $A$.

Let $Y$ be a generic hyperplane in $\C^{n+m}$ tangent to the manifold $\mbox{sing} A$ at some its point $s$. Let $B$ be a small ball centered at $s$, and $\tilde Y \in {\mathcal Reg}^{\C}$ be a hyperplane parallel to $Y$ and very close to it.

\begin{lemma}
\label{loclemma}
 Both groups participating in the map
\begin{equation}
\label{equat77}
 H_{n+m-1}(A \cap B, A\cap B \cap \tilde Y) \to H_{n+m-1}(A \cap B, (\AA \cup \tilde Y) \cap B)
\end{equation}
from the exact sequence of the triple $(A \cap B,  (\AA \cup \tilde Y) \cap B  ,  A \cap B \cap \tilde Y )$  are isomorphic to $\Z$, and this map is an isomorphism.
\label{lem17}
\end{lemma}

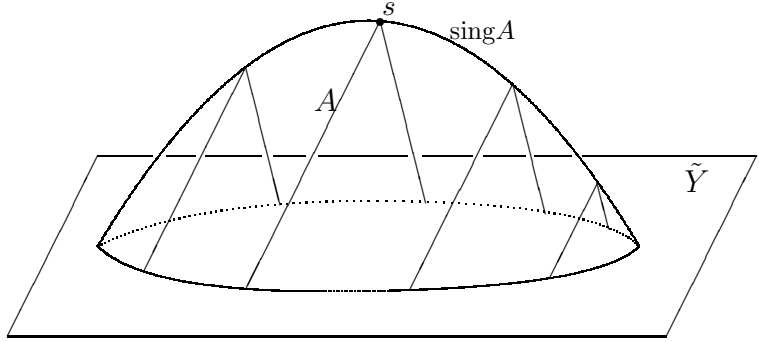
\begin{figure}
\begin{picture}(83,40)
\put(41.5,35.7){{\small $s$}}
\put(41.3,34.8){\circle*{1}}
\put(75,16){$\tilde Y$}
\put(49,33){{\footnotesize $\mbox{sing} A$}}
\put(34,25){$A$}
\bezier{700}(10,10)(40,60)(70,10)
\bezier{200}(10,10)(15,5)(35,5)
\bezier{200}(35,5)(65,5)(70,10)
\bezier{50}(10,10)(20,15)(40,15)
\bezier{50}(40,15)(65,15)(70,10)
\put(26.4,29.8){\line(-1,-2){11.33}}
\put(26.4,29.8){\line(1,-4){3.8}}
\put(41.3,35){\line(-1,-2){14.9}}
\put(41.3,35){\line(1,-4){5.05}}
\put(56,28){\line(-1,-2){11.38}}
\put(56,28){\line(1,-4){3.6}}
\put(65.35,17){\line(-1,-2){5.32}}
\put(65.35,17){\line(1,-4){1.26}}
\put(0,0){\line(1,0){73}}
\put(0,0){\line(1,2){10}}
\put(73,0){\line(1,2){10}}
\put(10,20){\line(1,0){6}}
\put(17.6,20){\line(1,0){3.2}}
\put(22.0,20){\line(1,0){6.2}}
\put(29.5,20){\line(1,0){3.5}}
\put(34.5,20){\line(1,0){9.8}}
\put(46.0,20){\line(1,0){5.5}}
\put(52.8,20){\line(1,0){4.6}}
\put(58.4,20){\line(1,0){4.1}}
\put(63.7,20){\line(1,0){19.3}}
\end{picture}
\caption{Fibration of $A$ close to the singular stratum}
\label{fi19}
\end{figure}

{\it Proof.} The space of choices of initial data participating in our lemma (consisting of the point $s \in \mbox{sing}A$, tangent plane $Y$, ball $B$, and neighbouring plane $\tilde Y \in {\mathcal Reg}^{\C}$) is irreducible, therefore it is enough to consider an arbitrary such collection of them. So we take for $s$ the point with coordinates $y_1=i\sqrt{1-\varepsilon^2}, y_2 = \dots = y_{m}=0$, choose the radius $\rho$ of the ball $B$ much smaller than $\varepsilon$,  and define $\tilde Y$ by the equation $y_1= i (\sqrt{1-\varepsilon^2}-\zeta)$ where $\zeta$ is a positive constant much smaller than $\rho$. We will show that both groups (\ref{equat77}) are then generated by one and the same relative cycle, namely by the part of cycle 2) from Lemma \ref{lem5} placed in the domain where \ $\mbox{Im\ } y_1 \in \left[\sqrt{1-\varepsilon^2}-\zeta , \sqrt{1-\varepsilon^2}\, \right]$.

The left-hand group in (\ref{equat77}) appears as follows.
The group $\tilde H_*(A \cap B)$ is trivial, as $A\cap B$ is homeomorphic to the cone over the point $s$.
The group $H_*(A \cap B \cap \tilde Y)$ can be reduced to the similar homology group in the transversal slice of $\mbox{sing} A$ by a kind of $(m-1)$-fold suspension, see \cite{GM}. Let us remind a realization of this reduction following \cite{APLT}, \S II.4 (see Fig. \ref{fi19}). Fiber the ball $B$ into its sections by a family of parallel $(n+1)$-dimensional complex planes transversal to $\mbox{sing} A$, along any of which the coordinates $y_2, \dots, y_m$ take some fixed values. Some of these planes are non-generic with respect to the hypersurface $A \cap\tilde Y$: this happens when the intersection point of such a plane and the variety $\mbox{sing} A$ belongs to $\tilde Y$. By the usual complex Morse lemma (applied to the restriction of the function $y_1$ to $\mbox{sing}A$), the variety $\mbox{sing} A \cap \tilde Y \cap B$ parameterizing the non-generic planes contains a sphere $S^{m-2}$ as a deformation retract. This sphere can be spanned by a $(m-1)$-dimensional disc $\Psi^{m-1},$ whose interior points lie in $\mbox{sing} A \cap B \setminus \tilde Y$; the class of this disc in
 $H_{m-1}(\mbox{\rm sing} A \cap B, \mbox{\rm sing} A \cap B \cap \tilde Y)$ is not equal to zero.

Let $\Upsilon$ be a generic complex $(n+1)$-dimensional plane from our family, which is transversal to $\mbox{sing} A$ at a point of $\Psi^{m-1}$, and $\Xi$ be an $i$-dimensional cycle in $A \cap B \cap \Upsilon \cap \tilde Y$. Then we can span an $(i+m-1)$-dimensional cycle in the variety $A \cap B \cap \tilde Y$ extending this cycle by the local triviality into all similar slices of this variety by the planes of this family intersecting $\mbox{sing} A$ in interior points of the disc $\Psi^{m-1}$ and contracting them  over the boundary points of this disc. This suspension operation defines an isomorphism $H_i(A \cap B \cap \Upsilon \cap \tilde Y) \to H_{i+m-1}(A \cap B \cap \tilde Y)$.

In our case, by Lemma \ref{lem4} the variety $ A \cap B \cap \Upsilon \cap \tilde Y$ is homeomorphic to a generic hyperplane section of the cone (\ref{quadrcone}), hence is homotopy equivalent to $S^{n-1}$. This sphere $S^{n-1}$ can be realized by the intersection of $A \cap \Upsilon \cap \tilde Y$ with cycle 2) from Lemma \ref{lem5}. Therefore $H_{n+m-2}(A \cap B \cap \tilde Y) \simeq \Z$, and  the equality of the first group in (\ref{equat77}) to $\Z$ follows by the exact sequence of the pair $(A \cap B, A \cap B \cap \tilde Y)$. By the construction a generator of this group is realized by the corresponding part of cycle 2) of Lemma \ref{lem5}.

Considering the second group (\ref{equat77}) we can replace the set $\AA$ by the complement in $A$ of a small closed tubular neighbourhood $T$ of the singular locus $\mbox{sing} A$; so we study the relative homology group \begin{equation}
\label{kunnet}
H_*(A \cap B,(A\setminus T) \cup (A \cap \tilde Y)) \cap B   ) \equiv H_*(A \cap T \cap B, A \cap (\partial T \cup \tilde Y) \cap B).\end{equation}

 By Lemma \ref{lem4} the variety $A \cap T \cap B$ is homeomorphic to the direct product of the $(2m-2)$-dimensional disc $\mbox{sing} A \cap B$ and the variety $A \cap D$ considered in Corollary \ref{cor2}. We can arrange the fibration of the tubular neighbourhood of $\mbox{sing}A$ in $B$ in such a way that each fiber either completely belongs to $\tilde Y$ or does not intersect it. Then our homeomorphism turns the pair $ (A \cap T \cap B, A \cap (\partial T \cup \tilde Y) \cap B)$ from (\ref{kunnet}) to the pair
$$((\mbox{sing}A \cap B) \times (A \cap D), ((\mbox{sing}A \cap \tilde Y \cap B) \times (A \cap D)) \cup ((\mbox{sing}A \cap B) \times (A \cap \partial D))).$$
The corresponding relative homology group (i.e. the right-hand group in (\ref{kunnet})) is therefore isomorphic to $\Z$ by Corollary \ref{cor2} and K\"unneth formula for the direct product of pairs $(\mbox{sing}A \cap B, \mbox{sing}A \cap \tilde Y \cap B)$ and $(A \cap D, A \cap \partial D)$. A generator of this group is again realized by cycle 2) from Lemma \ref{lem5} (intersected with the neighbourhood $T$).
\hfill $\Box$

\begin{lemma}
\label{lem6}
For any affine complex hyperplane $ X \in {\mathcal Reg}^{\C}$ in $\C^{n+m},$  we have isomorphisms
\begin{equation}H_{n+m}(\C^{n+m}, A \cup X) \simeq H_{n+m-1}(A \cup X) \simeq H_{n+m-1}(A,A \cap X)    \  .
\label{tit}
\end{equation}
All these groups are isomorphic to $\Z^6$.
\end{lemma}

\noindent
{\it Proof.}
The isomorphisms (\ref{tit}) follow from exact sequences of pairs.

By Thom isotopy theorem, all homology groups in (\ref{tit}) form coverings over the path-connected space ${\mathcal Reg}^{\C}$, so it is sufficient to calculate these groups for single hyperplane $X$ given by equation \ $x_1 +y_1=0$.
We will assume that $\varepsilon$ in (\ref{meq}) is small enough, in particular $\varepsilon < 1/\sqrt{2}$.

 The restriction of the function \ $x_1+ y_1$ \ to the stratified variety $A$ has six critical points. Four of them are Morse critical points placed in $\AA$, namely it are real points with $(x_1,y_1) = \pm (1,0) \pm (\varepsilon/\sqrt{2},\varepsilon/\sqrt{2})$, $x_2 = \dots =x_n=y_2= \dots = y_n=0$ (the signs $\pm$ are independent). In addition we have two critical points $\{x_1 = \dots =x_n=y_2= \dots =y_m=0, $ $ y_1 = \pm i \cdot \sqrt{1-\varepsilon^2}\}$ of the restriction of function $x_1+y_1$ to the singular locus (\ref{singloc}); they are Morse critical points of this function $x_1+y_1$ on the stratified variety $A$ in the  sense of \cite{GM}, \S 1.4.

\unitlength 0.90mm
\linethickness{0.4pt}
\begin{figure}
\begin{picture}(50,40)
\put(25,20){\circle{2}}
\put(0.6,20){\circle*{1.5}}
\put(49.5,20){\circle*{1.5}}
\put(5,20){\circle*{1.5}}
\put(45,20){\circle*{1.5}}
\put(25,0){\circle*{1.5}}
\put(25,40){\circle*{1.5}}
\put(24,20){\line(-1,0){18.5}}
\put(26,20){\line(1,0){18.5}}
\put(25,19){\line(0,-1){18.5}}
\put(25,21){\line(0,1){18.5}}
\bezier{150}(25.8,20.8)(37,30)(49.3,20.7)
\bezier{150}(24.2,19.2)(13,10)(0.7,19.3)
\end{picture}
\caption{Vanishing cycles for the model pair $(A,A \cap X)$}
\label{fi2}
\end{figure}
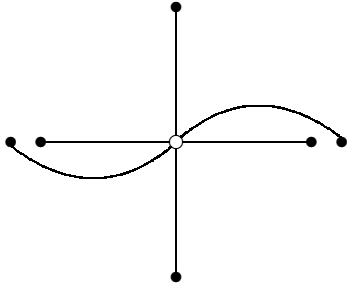

Let $B_j$, $j =1, \dots, 6$, be small balls in $\C^{n+m}$ around these six critical points. 
Consider a graph in $\C^1$ consisting of closed paths connecting the non-critical value 0 with all critical values  of the map $x_1+y_1: A \to \C^1$ as in Fig. \ref{fi2}. Denote these critical values by $z_j$, and denote by $\tilde z_j$ some points of this graph which are very close to its endpoints $z_j$; let $\tilde Y_j \in {\mathcal Reg}^{\C}$ be  six hyperplanes given by equations $x_1+y_1= \tilde z_j$.
 Denote by $\ddagger$ the pre-image in $A$ of this graph under the map $x_1+y_1$, and by $\breve \ddagger$ the pre-image of the same graph less the union of tiny segments $(\tilde z_j,z_j]$ at the ends of its legs.

The space $\ddagger$ is a deformation retract of $A$, and $A \cap X$ is a deformation retract of $\breve \ddagger$. This implies the first  isomorphism of the sequence
$$H_*(A,A \cap X) \simeq H_*(\ddagger, \breve \ddagger) \simeq 
\bigoplus_{j=1}^6 H_*(A \cap (x_1+y_1)^{-1}([\tilde z_j,z_j]), A \cap \tilde Y_j) \simeq $$
\begin{equation} \qquad \qquad \qquad
\simeq \bigoplus_{j=1}^6 H_{*}(A \cap B_j, A \cap \tilde Y_j \cap B_j) ;
\label{locliz}
\end{equation}
here the second equality is excision, and the third one follows from the local triviality of the restrictions of the map $x_1+y_1$ to the sets $(A \setminus B_j) \cap (x_1+y_1)^{-1}([\tilde z_j,z_j])$.

All six summands in (\ref{locliz}) are isomorphic to $\Z$ in dimension  $n+m-1$ and are trivial in all other dimensions.
For four summands related with non-vertical legs of the graph of Fig. \ref{fi2} this follows from the usual Morse lemma (and the corresponding generators of the group $H_{n+m-1}(A,A \cap X)$ are standard Lefschetz thimbles, see \cite{APLT}); for two vertical legs the same follows from Lemma \ref{loclemma}.    \hfill $\Box$

\medskip
\begin{corollary}
\label{coreppla}
For any $X \in {\mathcal Reg}^{\C}$, the group $H_{n+m-2}(A \cap X,\Q)$ is isomorphic to $\Q^4$.
\end{corollary}

\noindent{\it Proof.}
This follows immediately from lemmas \ref{lem5} and \ref{lem6} and exact sequence of the pair $(A,X)$. \hfill $\Box$

\medskip
Denote by $rH_*(A,A \cap X)$ the subgroup in $H_*(A,A \cap X)$ represented by relative cycles avoiding the singular locus of $A$, i.e. the image of the obvious map $H_*(\AA,\AA \cap X) \to H_*(A,A \cap X)$.

\begin{lemma} For any $X \in {\mathcal Reg}^{\C}$,  the group $rH_{n+m-1}(A,A \cap X)$ is isomorphic to $\Z^4$.
For $X$ considered in the proof of Lemma \ref{lem6} this group is generated by four Lefschetz thimbles defined by non-vertical paths in Fig. \ref{fi2}.
\label{lem8}
\end{lemma}

\noindent
{\it Proof.} The groups $rH_*(A,A\cap X)$ for all generic $X$ are isomorphic to one another, so let us take $X$ from the proof of Lemma \ref{lem6}.  By Lemma \ref{loclemma} (and the exact sequence of triple mentioned in its statement), the image of the group $H_{n+m-1}(\AA, \AA \cap X)$ in the group (\ref{locliz}) can be at most four-dimensional, since it does not contain non-zero elements of the summands corresponding to  two critical points from $\mbox{sing} A$. Four other generators of this sum can be realized by standard Lefschetz thimbles and hence belong to the image of this group. \hfill $\Box$

\section{Monodromy action}

Fundamental groups $\pi_1({\mathcal Reg}^{\C},X)$ and $\pi_1({\mathcal Reg}^{\C}_2,X)$   act by monodromy operators on all groups (\ref{tit}) and on the subgroup $rH_*(A,A \cap X)$ of the last of them; this action commutes with isomorphisms (\ref{tit}), cf. \cite{V14}. 

\begin{theorem}
\label{month}
The image of the group $\pi_1({\mathcal Reg}^{\C},X)$ in the group of automorphisms of the lattice $rH_*(A,A \cap X) \sim \Z^4$ under the monodromy representation is the Klein four-group $\Z_2 \oplus \Z_2$. For some concordant choice of orientations of basic Lefschetz thimbles generating this lattice, the automorphisms defined by this monodromy action preserve the set of these basic elements and act on it by permutations $\binom{1234}{2143}$ and $\binom{1234}{4321}$ and their compositions. The loops realised inside ${\mathcal Reg}_2^{\C}$ generate only the first of these permutations.
\end{theorem}

A proof of this theorem takes the rest of this section.
\medskip

Since the spaces ${\mathcal P}^{\C}$ and ${\mathcal P}^{\C}_2$ are simply-connected, these fundamental groups are generated by ``pinches'', i.e. the loops going first from the distinguished point $X$ to a smooth piece of the {\em discriminant set} $\Sigma \equiv {\mathcal P}^{\C} \setminus {\mathcal Reg}^{\C}$ of non-generic hyperplanes, then running a small circle around this piece, and coming back to $X$ along the first part of the loop.

The discriminant variety $\Sigma$ consists of four irreducible components $\Sigma_j$, $j=1,\dots, 4,$ formed respectively by hyperplanes

($\Sigma_1$) tangent to the variety $A$ at its non-singular finite points,

($\Sigma_2$) tangent to the variety $\mbox{sing} A$ at its finite points,

($\Sigma_3$)
defined by equations (\ref{eqhy}) with
\begin{equation}
\alpha_1^2 + \dots + \alpha_n^2+\gamma_1^2 + \dots + \gamma_m^2=0
\label{discr3}
\end{equation}
 (i.e. asymptotic to the smooth part of $A$: the closures in $\CP^{n+m}$ of these hyperplanes are tangent to the intersection of the closure of $A$ and the ``infinitely distant'' hyperplane $\CP^{n+m-1}_\infty \equiv \CP^{n+m} \setminus \C^{n+m}$), and

($\Sigma_4$) defined by equations (\ref{eqhy}) with
\begin{equation}
\gamma_1^2 + \dots + \gamma_m^2=0
\label{discr4}
\end{equation}
 (i.e. asymptotic with respect to  $\mbox{sing} A$).
\medskip

Let us study the action of pinches embracing these components  on the group $rH_{n+m-1}(A,A \cap X)$ .

\subsection{$\Sigma_1$: standard Picard--Lefschetz operator}
\begin{proposition} The pinches in ${\mathcal Reg}^{\C}$ going around smooth pieces of the component  $\Sigma_1$ act trivially on the group $rH_{n+m-1}(A,A \cap X)$.
\end{proposition}

\noindent
{\it Proof.}
Let $s$ be a regular point of $A$, at which the second fundamental form of $A$ is non-degenerate; let $Y \in \Sigma_1$ be the hyperplane tangent to $A$ at this point. We assume that $Y$ is not tangent to $A$ at any other points and does not belong to any other components of $\Sigma$ (which is true for almost all points $s \in A$). Then $\Sigma$ is smooth at the point $Y$. Let $L: \C^{n+m} \to \C$ be a linear (generally non-homogeneous) function such that $Y$ is defined by the equation $L=0$. Let $B$ be a small ball centered at $s$. Consider the one-parametric family of parallel hyperplanes $X_\lambda,$ $\lambda \in \C^1$, defined by equations $L(x,y) =\lambda$. For all $\lambda$ with sufficiently small $|\lambda| \neq 0$ these hyperplanes belong to ${\mathcal Reg}^{\C}$. Let $\lambda_0 \neq 0$ be such a value of the parameter $\lambda$ with very small $|\lambda_0|$. We need to calculate the monodromy action on $rH_{n+m-1}(A, A \cap X_{\lambda_0})$ defined by the circle in ${\mathcal Reg}^{\C}$ consisting of all hyperplanes $X_\lambda$ with $\lambda = e^{it} \lambda_0$, $t \in [0,2\pi]$. This action is described by the Picard--Lefschetz formula (\ref{plf}) formulated in the following terms.

 The group $H_{n+m-1}(A \cap B,A \cap B \cap X_{\lambda_0})$ is isomorphic to  $\Z$ and is generated by the {\em vanishing relative cycle} $\Delta(\lambda_0)$, see e.g. \cite{APLT}. The boundary map $H_{n+m-1}(A \cap B, A \cap B \cap X_{\lambda_0}) \to H_{n+m-2}(A \cap B \cap X_{\lambda_0})$ is an isomorphism, and the latter group is generated by the absolute cycle $\partial \Delta(\lambda)$.

Further, let $\Delta$ be any element of the group $H_{n+m-1}(A,A\cap X_{\lambda_0})$, and $\partial \Delta$ be its boundary in $H_{n+m-2}(A \cap X_{\lambda_0})$. Then the monodromy operator in question sends our class $\Delta$ to
\begin{equation}
\Delta + (-1)^{(n+m)(n+m-1)/2} \langle \partial \Delta, \partial \Delta(\lambda_0) \rangle \Delta(\lambda_0) \ ,
\label{plf}
\end{equation}
where $\langle \cdot, \cdot \rangle$ is the intersection index in $\AA \cap X_{\lambda_0}$. Therefore it remains to prove the following lemma.

\begin{lemma}
\label{lemseven}
 For $X \in {\mathcal Reg}^{\C}$,
any two elements of $rH_{n+m-1}(A, A \cap X)$ can be represented by relative cycles  $\Delta, \Delta' \subset \AA$ such that  $\langle \partial \Delta, \partial \Delta' \rangle =0$.
\end{lemma}

\unitlength 1.50mm
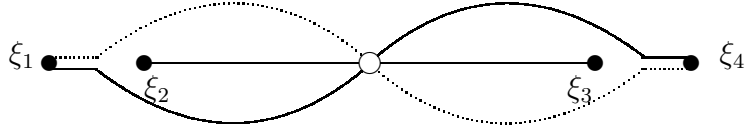
\begin{figure}
\begin{picture}(60,20)
\put(29,10){\circle{1.8}}
\put(0.6,10){\circle*{1.3}}
\put(57.5,10){\circle*{1.3}}
\put(9,10){\circle*{1.3}}
\put(49,10){\circle*{1.3}}
\put(28,10){\line(-1,0){18.5}}
\put(30,10){\line(1,0){18.5}}
\put(0.8,9.5){\line(1,0){4.1}}
\put(57.1,10.5){\line(-1,0){4}}
\bezier{800}(29.7,10.7)(41,20)(53.3,10.6)
\bezier{800}(28.3,9.3)(17,0)(4.7,9.5)
\bezier{55}(29.8,9.2)(41,00)(53.3,9.5)
\bezier{55}(28.2,10.8)(17,20)(4.7,10.5)
\bezier{10}(0.8,10.5)(2.8,10.5)(4.9,10.5)
\bezier{10}(57.1,9.5)(55.1,9.5)(53.1,9.5)
\put(-3,10){$\xi_1$}
\put(9,7){$\xi_2$}
\put(46.5,7){$\xi_3$}
\put(60,10){$\xi_4$}
\end{picture}
\caption{Vanishing cycles and their complex conjugates}
\label{fi4}
\end{figure}

\noindent
{\it Proof.} It is sufficient to prove this for the hyperplane $X$ used in the proof of Lemma \ref{lem6}, and for basic Lefschetz thimbles issuing from the four critical points of the function $x_1+y_1: \AA \to \C$, see Figs. \ref{fi2}, \ref{fi4}. All these critical points are real, and  the Morse indices of their restrictions to the real part $\partial W$ of $A$ are even.
Let $\xi_j$, $j\in \{1,\dots, 4\}$, be critical values $\pm 1 \pm \sqrt{2}\varepsilon$ of this restriction ordered by their increase. Consider some such value $\xi_j$, and let $\xi'$ be a non-critical real value of this restriction which is very close to $\xi_j$. Let $X_{\xi'}$ be the hyperplane in $\C^{n+m}$ defined by equation $x_1+y_1=\xi'$. The vanishing relative cycle in $(\AA,\AA \cap X_{\xi'})$ corresponding to the critical point with value $\xi_j$ is located in a small ball $B$ around this critical point and generates the group $H_{n+m-1}(\AA \cap B, \AA \cap X_\xi \cap B)$. This vanishing cycle $\Delta_j(\xi')$ can be chosen to be invariant under the complex conjugation; moreover, since the Morse index of the function $(x_1+y_1)|_{\partial W}$ is even, this conjugation preserves the orientation of this cycle and homology class of its boundary in $H_{n+m-2}(\AA \cap X_{\xi'})$. By isotopy considerations, this invariance holds also for all cycles $\Delta_j(\tilde \xi)$ in similar groups $H_{n+m-2}(\AA \cap X_{\tilde \xi}),$ where $\tilde \xi$ belong to the same interval of real non-critical values as $\xi'$, and $\Delta_j(\tilde \xi)$ are obtained from our vanishing cycle $\Delta_j(\xi')$ by the Gauss-Manin connection over this interval.

On the other hand, the complex conjugation changes the orientation of the $2(n+m-2)$-dimensional manifold $\AA \cap X_{\tilde \xi}$, therefore the intersection index of any two invariant cycles is opposite to itself, and hence is equal to 0. In particular, $\langle \partial \Delta_i(\xi),\partial \Delta_{i+1}(\xi) \rangle=0$, if $\xi$ is a non-critical value in the interval bounded by the critical values $\xi_i$ and $\xi_{i+1}$, and $ \Delta_j(\xi)$, $j=i$ or $i+1$, is the Lefschetz thimble in $(\AA,X_{\xi})$ defined by the segment connecting $\xi$ with the corresponding endpoint of this interval.

We obtain immediately the equality $\langle \partial \Delta_2, \partial \Delta_3 \rangle =0$, where $\Delta_2$, $\Delta_3$ are the Lefschetz thimbles in $(\AA,\AA \cap X)$ defined by segments in Fig. \ref{fi4} connecting 0 and critical values $\pm (1-\sqrt{2}\varepsilon)$.

Moreover, the intersection index does not change if we deform continuously in $\C^1$ the pattern consisting of the non-critical value $\xi'$ and two paths connecting it to the neighbouring (fixed) critical values in such a way that all this construction has no additional meetings with the set of critical values during this movement. For instance, the union of two solid paths in Fig. \ref{fi4} connecting the critical values $\xi_1$ and $ \xi_2$ to the non-critical value 0 can be obtained by such a deformation from two parts of the interval $(\xi_1,\xi_2) \in \R^1$ connecting the same critical values to some point of this interval.
This implies equality  $\langle \partial \Delta_1, \partial \Delta_2 \rangle = 0$ and, analogously, $\langle \partial \Delta_3, \partial \Delta_4 \rangle =0$.

Further, denote by $\tilde \Delta_1$ and $\tilde \Delta_4$ Lefschetz thimbles in  $(\AA,\AA \cap X)$ which are complex conjugate to $\Delta_1$ and $\Delta_4$ respectively, and therefore can be defined by the pointed paths in Fig. \ref{fi4}.
The cycles $\partial \Delta_1 + \partial \tilde \Delta_1$ and  $\partial \Delta_4 + \partial \tilde \Delta_4$ are invariant under the complex conjugation, therefore we have $$\langle \partial \Delta_1 + \partial \tilde \Delta_1, \partial \Delta_3 \rangle =0, \quad
\langle \partial \Delta_2, \partial \Delta_4 + \partial \tilde \Delta_4 \rangle =0,$$ $$ \langle \partial \Delta_1 + \partial \tilde \Delta_1,  \partial \Delta_4 + \partial \tilde \Delta_4 \rangle =0.$$
By Picard--Lefschetz formula $\partial \tilde \Delta_1$ is equal to $\partial \Delta_1 \pm \langle \partial \Delta_1, \partial \Delta_2 \rangle \partial \Delta_2$, and hence (by the equality  $\langle \partial \Delta_1, \partial \Delta_2 \rangle =0$  proved in the previous paragraph) to $\partial \Delta_1$.   In a similar way, $\partial \tilde \Delta_4 = \partial \Delta_4$. Therefore the previous three equalities reduce to  $\langle \partial \Delta_1, \partial \Delta_3 \rangle =0,$
$\langle \partial \Delta_2, \partial \Delta_4 \rangle =0,$ and $\langle \partial \Delta_1 ,  \partial \Delta_4 \rangle =0.$

Finally, all self-intersection indices $\langle \partial \Delta_j, \partial \Delta_j \rangle$ of odd-dimensional cycles are trivial.
\hfill $\Box$ $\Box$

\subsection{$\Sigma_2$: tangents to the singular locus.}

Let $Y$ be a generic hyperplane in $\C^{n+m}$ tangent to $\mbox{sing} A$ at some point $s$, so that its linear (generally, non-homogeneous) equation $L=0$ defines a Morse function on the stratified variety $A$ in  the sense of \cite{GM}, and $Y$ does not belong to other local branches of $\Sigma$. Let $Y_\xi$, $\xi \in \C$, be the family of hyperplanes parallel to $Y \equiv Y_0$, which are defined by the equations $L(x,y)=\xi$. Let $B$ be a small ball in $\C^{n+m}$ centered at $s$. Let $\xi_0\neq 0$ be a number such that $|\xi_0|$ is much smaller than the radius of the ball $B$. We consider the monodromy operator acting on the group
\begin{equation}H_{n+m-1}(A, A \cap Y_{\xi_0})
\label{hc3}
\end{equation}
and defined by the Gauss-Manin connection over the circle in ${\mathcal Reg}^{\C}$ consisting of hyperplanes $Y_\xi$ where $\xi$ runs the circle $e^{it} \xi_0$, $t \in [0,2\pi]$.

\begin{proposition}
This monodromy operator acts trivially on the subgroup $rH_{n+m-1}(A,A \cap Y_{\xi_0})$.
\end{proposition}

\noindent
{\it Proof.} The family of varieties $(A \cap Y_\xi \setminus B)$ forms a trivial fiber bundle over the disc in $\C$ consisting of values $\xi$ with $|\xi| \leq |\xi_0|$, therefore
 this monodromy operator adds to any element $\Delta$ of the group (\ref{hc3}) some element (depending linearly on $\Delta$) which can be realized by a chain inside $B$. Moreover, if $\Delta \in rH_{n+m-1}(A,A \cap Y_{\xi_0})$ then this chain can be realized inside $\AA \cap B$, and hence its homology class in $H_{n+m-1}(A \cap B, A \cap B \cap Y_{\xi_0})$ belongs to the image of the homomorphism
\begin{equation}
\label{equat7}
H_{n+m-1}(\AA \cap B, \AA \cap B \cap Y_{\xi_0}) \to H_{n+m-1}(A \cap B, A\cap B \cap Y_{\xi_0})
\end{equation}
of the exact sequence of the triple $(A\cap B, (\AA \cup Y_{\xi_0}) \cap B, A \cap B \cap Y_{\xi_0})$. This image is trivial by Lemma \ref{loclemma}. \hfill $\Box$

\subsection{$\Sigma_3$ and $\Sigma_4$: asymptotic hyperplanes}
\label{asy}

The space ${\mathcal P}^{\C}$ augmented with the point corresponding to the ``infinitely distant'' hyperplane in $\CP^{n+m}$ is itself isomorphic to the $(n+m)$-dimensional complex projective space. Any generic two-dimensional projective subspace in it intersects the varieties $\Sigma_3$ and $\Sigma_4$
(distinguished by conditions (\ref{discr3}) and (\ref{discr4})) along two degree 2 complex curves in general position, in particular these curves have exactly four transversal intersection points. (If $m>2$ then both these curves are non-singular, and in the case $m=2$ the second of them splits into two lines). By Zariski theorem, the fundamental group of the complement of the union of these curves in this 2-subspace (and hence also of the complement of  variety $\Sigma_3 \cup \Sigma_4$ in entire ${\mathcal P}^{\C}$) is isomorphic to $\Z^2$, in particular is commutative. Therefore (and since components $\Sigma_1$ and $\Sigma_2$ do not contribute to the monodromy action, as is shown in two previous subsections) we can calculate independently the action on $rH_{n+m-1}(A,A \cap X)$ of loops embracing these components, and not take care on the choice of the distinguished point in ${\mathcal Reg}^{\C}$. We will do it for two loops, whose linking numbers with varieties $\Sigma_3$ and $\Sigma_4$ are equal to $(1,0)$ and $(1,1)$.
\smallskip

{\it The first of these loops}  will be realized inside the subspace ${\mathcal P}^{\C}_2 \subset {\mathcal P}^{\C}$ (see \S \ref{mlmm}). 
This subspace meets the discriminant variety $\Sigma_4$ at the points of its singular part only; the complement of this variety in this subspace is simply-connected, so that all loops in ${\mathcal Reg}_2^{\C}$ have zero linking number with $\Sigma_4$. They generate a subgroup isomorphic to $\Z$ in the lattice $\pi_1({\mathcal P}^{\C} \setminus (\Sigma_3 \cup \Sigma_4)) \simeq \Z^2$, and the linking numbers with $\Sigma_3$ separate the points of this subgroup.

This subspace intersects the set $\Sigma_3$ along two lines consisting of hyperplanes given by equations $x_1=\pm i y_1 + \beta$ with arbitrary $\beta$. The space ${\mathcal P}_2^{\C} \setminus \Sigma_3$ can be projected to $\CP^1\setminus \{i,-i\}$ by sending any plane with equation (\ref{hhh}) to the number $p/q$. A generator of the fundamental group of this space is provided by any loop along which these numbers run the real axis $\RP^1 \subset \CP^1$, or by any loop sufficiently close to it.

\begin{figure}
\begin{picture}(41,40)
\bezier{250}(7,20)(7,33)(20,33)
\bezier{250}(20,33)(33,33)(33,20)
\bezier{250}(33,20)(33,7)(20,7)
\bezier{250}(20,7)(7,7)(7,20)
\put(29.8,29.8){\circle*{1}}
\put(10.2,10.2){\circle*{1}}
\put(40,18.5){\line(-1,1){21.5}}
\put(0,21.5){\line(1,-1){21.5}}
\put(29.8,29.8){\circle{10}}
\put(10.2,10.2){\circle{10}}
\put(20,20){\circle*{1.2}}
\put(17.3,17){{\small (1,0)}}
\bezier{70}(13,33)(8.5,31.5)(7,27)
\put(7,27){\vector(-1,-3){0.2}}
\put(22,-1){\small $X(5\pi/4)$}
\put(38,16){\small $X(\pi/4)$}
\put(29.5,8.5){$A$}
\put(20,20){\vector(1,0){13}}
\put(30,20.6){$\varepsilon$}
\end{picture}
\caption{A loop of third type}
\label{pl3}
\end{figure}
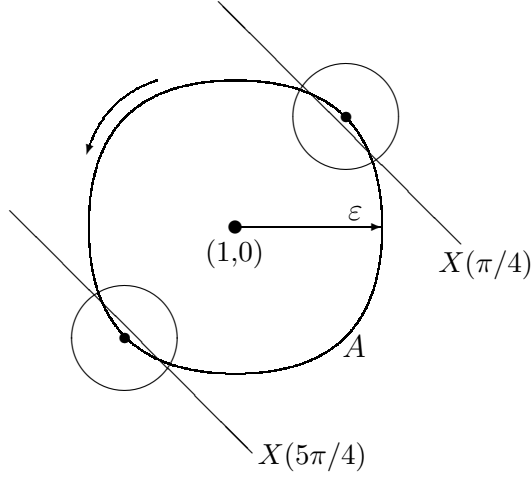

To realize such a loop, let us draw a path on the surface $\partial W$, connecting two critical points of the restriction of function
$x_1+y_1$ to $A$ inside the plane $\R^2$ defined by conditions
\begin{equation}
 x_2= \dots=x_n=y_2= \dots = y_m=0,
\label{r2}
\end{equation}
see Fig. \ref{pl3}: the coordinates $x_1$ and $y_1$ of the points of this path move as $x_1(\tau) = 1 + \varepsilon \cos \tau,$ $y_1(\tau) = \varepsilon \sin \tau$, $\tau \in [\pi/4, 5\pi/4]$. The tangent planes of $A$ at all these points can be defined by equations with real coefficients, in particular are distant from the set $\Sigma_3$. Let us deform this path slightly inside $A \cap \C^2$ (where $\C^2$ is the complexification of the plane $\R^2$ distinguished by (\ref{r2})) in such a way that the tangent hyperplanes of $A$ at all points $\eta(\tau),$ $\tau \in [\pi/4, 5\pi/4]$, of the obtained path are non-singular points of $\Sigma$; in particular all these points $\eta(\tau)$ are {\em non-parabolic} points of $A$, i.e. the second fundamental form of $A$ is non-degenerate at all of them. 

(There are exactly three such points of the initial real path in $\partial W$ where we need to move slightly to the complex domain. One of them corresponds to $\tau=\pi/2$, where the second fundamental form of $\partial W$ degenerates, the second to $\tau=\pi/2+ \arcsin \varepsilon$, where an extra tangency with $A$ at a distant point occurs, and the third one to $\tau=\pi$, where the tangent plane belongs to $\Sigma_4$.)

After that, let us pave a path in ${\mathcal Reg}^{\C}_2$ consisting of hyperplanes $X(\tau)$ which are parallel to the tangent hyperplanes of $A$ at the corresponding points $\eta(\tau)$ of our path, and are extremely close to these tangent hyperplanes.
For any point $\eta(\tau) \in A \cap \C^2$ of the first path and a small ball $B_\tau$ around this point, we have $H_{n+m-1}(A \cap B_\tau, A \cap B_\tau \cap X(\tau)) \simeq \Z$. Denote by $ \delta(\tau)$ the vanishing relative cycles generating these groups and supplied with some orientations depending continuously on $\tau$.

The Gauss--Manin connection over our path moves these vanishing cycles ones into the others, in particular it moves the vanishing cycle $\delta(\pi/4)$ related with the starting point to the vanishing cycle $\delta(5\pi/4)$ related with the final one.

Consider two paths in ${\mathcal Reg}_2^{\C}$ connecting the hyperplanes
$X(\pi/4)$ and $X(5\pi/4)$ with the hyperplane $X_0 = \{x_1+y_1=0\}$ and consisting of hyperplanes defined by equations $x_1+y_1=\lambda$, where $\lambda$ runs the solid paths in the right-hand part of Fig. \ref{fi4}. Deforming our two vanishing  relative cycles $\delta(\pi/4)$ and $\delta(5\pi/4)$ by the local triviality of the fiber bundle of pairs $(A, A \cap X)$
over these paths, we obtain exactly the Lefschetz thimbles  $\Delta_4$ and $\Delta_3$ respectively. Therefore, the loop in ${\mathcal Reg}_2^{\C}$ consisting of these two paths (passed in appropriate directions) and our path $\{X(\tau)\},$ $\tau \in [\pi/4,5\pi/4],$ moves the class $\Delta_4 \in rH_{n+m-1}(A,A \cap X_0)$ into $\Delta_3$.

 This loop is a generator of the group $\pi_1({\mathcal P}^{\C} \setminus \Sigma_3, X_0)$. A different realization of the same generator is provided by the union of the same two paths (passed in opposite directions) and the family of escorting hyperplanes $X(\tau)$ similar to the one considered above but corresponding to $\tau \in [5\pi/4,9\pi/4]$. By the same considerations as above the latter family moves the vanishing cycle $\delta(5\pi/4)$ to $\delta(\pi/4)$ or $-\delta(\pi/4)$, and hence the entire loop 
moves the class $\Delta_3$ to $\Delta_4$ (or, respectively, to $-\Delta_4$).

 It remains to prove that the correct choice of the sign is $+$. This property is equivalent to the assertion that the transportation of vanishing cycles $\delta(\cdot)$ along entire circle, composed of two paths $\{\eta(\tau)\}$, $\tau \in [\pi/4,5\pi/4]$ and $\tau \in [5\pi/4,9\pi/4]$, moves the homology class of the vanishing cycle $\delta(\tau)$ to itself and not to minus itself. Proving this we can start from an arbitrary point of the circle. 

To do it, notice that we may simplify our paths $\{\eta(\tau)\}$ and $\{X(\tau)\}$, not moving into the complex domain at the points $\tau = \pi$ and $2\pi$: indeed, in these cases the degeneration of the topological type of the pair $(A, A \cap X(\tau))$ happens far away from the support of our cycles $\delta(\tau)$, and the transportation of the cycles $\delta(\tau)$ along these straightened paths  gives the same result. Further, we can realize this family of cycles $\delta(\tau)$, $\tau \in [0,2\pi]$,   in such a way that for any two points $\tau, \tau' $ with $\tau + \tau' = 2\pi$ the corresponding cycles $\delta(\tau)$ and $\delta(\tau')$ will be symmetric to one another with respect to the involution in $\C^{n+m}$ multiplying the coordinate $y_1$ by $-1$. Indeed, we can choose both cycles corresponding to $\tau=\pi$ and $\tau=0 \equiv 2\pi$ to be symmetric to themselves, then realize arbitrarily  the family of cycles $\delta(\tau)$ connecting them over $\tau \in [0,\pi]$, and realize the family for $\tau \in [\pi,2\pi]$ by cycles symmetric to the corresponding cycles of the previous family. Then the transportations of the cycle $\delta(\pi)$ over the half-circles $\tau \in [0,\pi]$ and $[\pi,2\pi]$ give the same results, which is equivalent to our assertion.

So, our generator of the group $\pi_1({\mathcal P}^{\C} \setminus \Sigma_3, X_0)$ permutes basic classes $\Delta_3$ and $\Delta_4$ of the group $rH_{n+m-1}(A,A \cap X_0)$. In absolutely the same way, it permutes the classes $\Delta_1$ and $\Delta_2$; so it realizes the first permutation indicated in Theorem \ref{month}.
\medskip

{\it The second loop} in ${\mathcal Reg}^{\C}$
is provided by the family of hyperplanes  defined by real equations 
$$ (x_1+y_1) \cos \gamma + (x_2+ y_2) \sin \gamma =0 \ , \quad \gamma \in [0,\pi] \ . $$ It continuously permutes four critical points of the restrictions of functions $\pm(x_1+y_1)$ to the hypersurface $\partial W$, and also permutes the corresponding Lefschetz thimbles $\Delta_2 \leftrightarrow \Delta_3$ and $\Delta_1 \leftrightarrow \Delta_4$. The linking numbers of this loop with varieties (\ref{discr3}) and (\ref{discr4}) are equal to 1. \hfill $\Box$

\section{Lefschetz thimbles and  integration contours}

\begin{lemma}
\label{onlyfour}
If the restriction of a linear function $\R^{n+m} \to \R$ to $\partial W$ is Morse, then it has exactly four critical points.
\end{lemma}

Indeed, the number of its critical points cannot be smaller than four by the Morse inequality, and cannot be bigger because even the restriction of the complexification of a generic linear function to $\AA$ has no more than four critical points. \hfill $\Box$ \medskip

By Lemma \ref{lem8}, for any $X \in {\mathcal Reg}^{\C}$ the subgroup $rH_{n+m-1}(A, A \cap X) \subset H_{n+m-1}(A, A \cap X)$ is isomorphic to $\Z^4$ and is generated by the classes $\Delta_j$ of four Lefschetz thimbles associated with the critical points of the restriction to the smooth part of $A$ of the linear function $L$ vanishing on the hyperplane $X$. 

By Lemma \ref{lemseven}, the homology classes of these thimbles (defined up to a choice of their orientations) do not depend on  the paths connecting the critical values of $L|_{\AA}$ with $0$. Also, we can choose their orientations in concordant way so that they constitute one orbit of the group $\pi_1({\mathcal Reg}^{\C},X)$. (This choice of orientations is unique up to a simultaneous change of all of them).

Denote by $\Xi_j \in H_{n+m}(\C^{n+m},A \cup X)$ the preimages of these basic classes $\Delta_j$ under the composite isomorphism (\ref{tit}). 

Suppose now that the hyperplane $X$ is real, $X \in {\mathcal Reg}$, and the  linear function $L:\R^{n+m} \to \R$  vanishing on $X$ is Morse on $\partial W$. Let us fix some orientation of the space $\R^{n+m}$.

\begin{proposition}[see e.g. \cite{V14}]
\label{old}
The homology class in $H_{n+m}(\C^{n+m}, A \cup X)$ of the cycle  $W \cap L^{-1}((-\infty,0])$ $($respectively, $W \cap L^{-1}([0,+\infty)))$ oriented by the chosen orientation of \ $\R^{n+m}$ is equal to the sum of these cycles $\Xi_j$, defined by the paths connecting $0$ with all negative  $($respectively, positive$)$ critical values and  taken with appropriate coefficients equal to 1 or $-1$. \hfill $\Box$
\end{proposition}

\begin{proposition}
For some concordant choice of orientations of cycles $\Delta_j$, all four coefficients, with which the cycles $\Xi_j$ participate in the two sums mentioned in Proposition \ref{old}, are equal to $+1$.
\end{proposition}

Indeed, the sum of these two sums is equal to the class of the body $\partial W$, and hence does not depend on the hyperplane $X$. Therefore it is invariant under the action of the group $\pi_1({\mathcal Reg}^{\C},X)$. This is possible only if all coefficients are equal to one another. We can make all of them to be equal to $+1$ by the choice of one of two sets of concordant orientations.  \hfill $\Box$ 
\medskip

Now define the four-valued analytic function $\Psi$ on the space ${\mathcal P}^{\C}$, whose four values at any point $X \in {\mathcal Reg}^{\C}$ are equal to the integrals of the form (\ref{volform})  along all cycles $\Xi_j$. This function satisfies the conditions of Theorem \ref{redmt}. \hfill $\Box$

\end{document}